\documentclass[english]{article}

\usepackage[OT1]{fontenc}
\usepackage[latin9]{inputenc}
\usepackage{amsmath}
\usepackage{amsthm}

\makeatletter
\theoremstyle{plain}
\newtheorem{thm}{\protect\theoremname}
\theoremstyle{remark}
\newtheorem{rem}[thm]{\protect\remarkname}
\theoremstyle{plain}
\newtheorem{cor}[thm]{\protect\corollaryname}
\theoremstyle{plain}
\newtheorem{lem}[thm]{\protect\lemmaname}
\ifx\proof\undefined
\newenvironment{proof}[1][\protect\proofname]{\par
	\normalfont\topsep6\p@\@plus6\p@\relax
	\trivlist
	\itemindent\parindent
	\item[\hskip\labelsep\scshape #1]\ignorespaces
}{%
	\endtrivlist\@endpefalse
}
\providecommand{\proofname}{Proof}
\fi
\theoremstyle{definition}
\newtheorem{example}[thm]{\protect\examplename}
\newcommand{\lyxaddress}[1]{
	\par {\raggedright #1
	\vspace{1.4em}
	\noindent\par}
}

\@ifundefined{date}{}{\date{}}
\makeatother

\usepackage{babel}
\providecommand{\corollaryname}{Corollary}
\providecommand{\examplename}{Example}
\providecommand{\lemmaname}{Lemma}
\providecommand{\remarkname}{Remark}
\providecommand{\theoremname}{Theorem}

\begin{document}
\title{Gap theorems in Yang-Mills theory for complete four-dimensional manifolds
with a weighted Poincaré inequality}
\author{Matheus Vieira}
\maketitle
\begin{abstract}
In this paper we prove  $L^{\infty}$ type gap theorems in Yang-Mills
theory for complete four-dimensional manifolds with a weighted Poincaré
inequality. We apply the theorems to a broad class of complete manifolds
satisfying weighted Poincaré inequalities. In particular, we obtain
a gap theorem on the Euclidean space without assuming finite Yang-Mills
energy. We also prove an $L^{\infty}$ characterization of the BPST
instanton centered at the origin with unit scale.
\end{abstract}

\section{Introduction}

Suppose we have a Riemannian vector bundle with structure group $G\subset O\left(N\right)$
on a Riemannian manifold $X$ and a metric connection $A$ on the
bundle with curvature $F$. The connection $A$ is called a Yang-Mills
connection if the divergence of the curvature $F$ is zero, which
is the Euler-Lagrange equation of the Yang-Mills functional $YM\left(A\right)=\int_{X}\left|F\right|^{2}$.
When the manifold is four-dimensional, the connection $A$ is called
an instanton if the self-dual curvature $F^{+}$ is zero or the anti-self-dual
curvature $F^{-}$ is zero, where $F^{\pm}=1/2\left(F\pm*F\right)$.
Any instanton is a Yang-Mills connection. In \cite{AHDM1978} and
\cite{BPST1975} Atiyah-Drinfeld-Hitchin-Manin and Belavin-Polyakov-Schwartz-Tyupkin
found an explicit $SU\left(2\right)$ instanton, known as the BPST
 instanton. On the other hand, there are Yang-Mills connections that
are not instantons. In \cite{SSU1989} Sibner-Sibner-Uhlenbeck proved
the existence of $SU\left(2\right)$ Yang-Mills connections that are
not instantons (see also Bor \cite{B1992}, Parker \cite{P1992} and
Sadun-Segert \cite{SS1992}). For larger structure groups such as
$SU\left(4\right)$ we can easily find Yang-Mills connections that
are not instantons (see for example Section 2 in \cite{VN2008}).

By a gap theorem for the self-dual curvature $F^{+}$ of a Yang-Mills
connection we mean the following: under certain assumptions on the
four-dimensional manifold, if the self-dual curvature $F^{+}$ of
a Yang-Mills connection has a sufficiently small norm (for example
the $L^{2}$ norm or the $L^{\infty}$ norm), then the self-dual curvature
$F^{+}$ is zero. In \cite{BLS1979} and \cite{BL1981} Bourguignon-Lawson-Simons
proved an $L^{\infty}$ gap theorem for the self-dual curvature $F^{+}$
of a Yang-Mills connection on compact four-dimensional manifolds with
a certain positive curvature (for the extension to complete manifolds
see Shen \cite{S1982}). In \cite{M1982} MinOo proved an $L^{2}$
gap theorem for the self-dual curvature $F^{+}$ of a Yang-Mills connection
on these manifolds (see also Parker \cite{P1982}, and for the extension
to complete manifolds see Dodziuk-MinOo \cite{DM1982}). In \cite{F2016}
Feehan proved an $L^{2}$ gap theorem for the self-dual curvature
$F^{+}$ of a Yang-Mills connection on compact four-dimensional manifolds
with a good metric. In \cite{GKS2018} Gursky-Kelleher-Streets proved
an $L^{2}$ gap theorem for the self-dual curvature $F^{+}$ of a
Yang-Mills connection on compact four-dimensional manifolds with positive
Yamabe invariant. A result of this type for general compact four-dimensional
manifolds is still an open problem. In \cite{V2024} we extended the
gap theorem of Gursky-Kelleher-Streets to complete manifolds and we
also described the equality in the gap theorem in terms of the BPST
instanton, which is interesting even for compact manifolds. By a gap
theorem for the full curvature $F$ of a Yang-Mills connection we
mean the following: under certain assumptions on the $n$-dimensional
manifold, if the curvature $F$ of a Yang-Mills connection has a sufficiently
small norm (for example the $L^{n/2}$ norm or $L^{\infty}$ norm),
then the curvature $F$ is zero. See for example Price \cite{P1983},
Gerhardt \cite{G2010}, Zhou \cite{Z2015} and \cite{Z2016}, Feehan
\cite{F2017} and \cite{F2024}, and the above references. When the
manifold is four-dimensional, gap theorems for the self-dual curvature
$F^{+}$ and the anti-self-dual curvature $F^{-}$ are more important
than gap theorems for the full curvature $F$.

In this paper we extend the gap theorem of Bourguignon-Lawson-Simons
to complete four-dimensional manifolds with a weighted Poincaré inequality.
Using the weighted Poincaré inequality we can overcome the assumption
of positive curvature in the work of Bourguignon-Lawson-Simons. In
this way we cover space forms, complex space forms, manifolds with
positive spectrum, stable minimal hypersurfaces, and several other
geometric settings. Unlike previous works, in Theorem \ref{thm:vol}
we do not assume finite Yang-Mills energy. In Corollary \ref{cor:app}
we obtain $L^{\infty}$ gap theorems for all constant curvature models
(real and complex). In particular, we prove that a Yang-Mills connection
on $R^{4}$ satisfying $\left|r^{2}F^{+}\right|_{L^{\infty}}\leq2/\gamma$
must be anti-self-dual, without assuming finite Yang-Mills energy.
Another contribution related to this corollary is that we find a weighted
Poincaré inequality that improves the bottom of the spectrum of the
Laplacian of the complex hyperbolic space (Lemma \ref{lem:wpichn}).
Finally, in Theorem \ref{thm:adhm} we give an $L^{\infty}$ characterization
of the BPST instanton centered at the origin with unit scale.

Suppose we have an $L_{loc}^{1}$ function $q$ on a complete Riemannian
manifold $X$. The manifold $X$ is said to satisfy a weighted Poincaré
inequality with a weight function $q$ if the following inequality
is valid:
\begin{equation}
\int_{X}q\phi^{2}\leq\int_{X}\left|\nabla\phi\right|^{2}\label{eq:wpi}
\end{equation}
for all smooth functions $\phi$ on the manifold $X$ with compact
support. In \cite{C1997} and \cite{LW2006} Carron and Li-Wang applied
weighted Poincaré inequalities to Riemannian manifolds and submanifolds.
In Section 2 we give many examples of manifolds with a weighted Poincaré
inequality (the Euclidean space, the hyperbolic space, the complex
hyperbolic space, manifolds with positive spectrum, stable minimal
hypersurfaces in a general Riemannian manifold and minimal hypersurfaces
in the Euclidean space) and we also give  more references on this
subject.

In this paper we prove two $L^{\infty}$ type gap theorems for the
self-dual curvature $F^{+}$ of a Yang-Mills connection on complete
four-dimensional manifolds with a weighted Poincaré inequality. The
theorems have different assumptions and complementary applications.

We denote the distance function by $r\left(\cdot\right)=dist\left(\cdot,x_{0}\right)$.
We denote the open ball by $B_{R}=\left\{ r<R\right\} $. We denote
by $\lambda_{\max}\left(W^{+}\right)$ the largest eigenvalue of the
self-dual Weyl curvature $W^{+}=1/2\left(W+*W\right)$. We denote
by $\gamma$ a constant depending only on the structure group $G$
(see Lemma \ref{lem:bochner}).

In the next result, assuming that the manifold has volume growth $vol\left(B_{R}\right)=O\left(R^{\alpha}\right)$
and the self-dual curvature $F^{+}$ has growth $\left|F^{+}\right|=O\left(r^{2-\alpha}\right)$,
we prove an $L^{\infty}$ type gap theorem for the self-dual curvature
$F^{+}$ without any assumption on its integral. To our knowledge,
this is the first gap theorem for complete manifolds that does not
impose any integral assumption on the curvature.
\begin{thm}
\label{thm:vol}Suppose a complete four-dimensional Riemannian manifold
$X$, with scalar curvature $S$ and Weyl curvature $W$, satisfies
the weighted Poincaré inequality (\ref{eq:wpi}) with a weight function
$q$ (smooth on $X\setminus\left\{ x_{0}\right\} $) and has volume
growth $vol\left(B_{R}\right)=O\left(R^{\alpha}\right)$ for some
constant $\alpha>0$. Given a Yang-Mills connection $A$ with curvature
$F$ and structure group $G\subset O\left(N\right)$ ($N\geq3$) on
the manifold $X$, suppose the self-dual curvature $F^{+}$ has growth
$\left|F^{+}\right|=O\left(r^{2-\alpha}\right)$ as $r\to\infty$
and satisfies the pointwise inequality
\begin{equation}
\left|F^{+}\right|\leq\left(1/\gamma\right)\left(2q+\left(1/3\right)S-2\lambda_{\max}\left(W^{+}\right)\right)\,\,\,\,\,on\,\,\,X\setminus\left\{ x_{0}\right\} .\label{eq:volineq}
\end{equation}
Then either inequality (\ref{eq:volineq}) is an equality (on $X\setminus\left\{ x_{0}\right\} $)
or the self-dual curvature $F^{+}$ is identically zero.
\end{thm}
In the theorem above we also assume that the right hand side of inequality
(\ref{eq:volineq}) is nonnegative and not identically zero. We remark
that the growth assumption $\left|F^{+}\right|=O\left(r^{2-\alpha}\right)$
is what guarantees that the estimate in the proof involving the logarithmic
cutoff goes to zero.

In the next result, assuming that the weight function $q$ has growth
$q=O\left(r^{2}\right)$ and the self-dual curvature $F^{+}$ is in
the space $L^{2p}$, we prove an $L^{\infty}$ type gap theorem for
the self-dual curvature $F^{+}$ without any assumption on the volume
of the manifold.
\begin{thm}
\label{thm:int}Suppose a complete four-dimensional Riemannian manifold
$X$, with scalar curvature $S$ and Weyl curvature $W$, satisfies
the weighted Poincaré inequality (\ref{eq:wpi}) with a weight function
$q$ (smooth on $X\setminus\left\{ x_{0}\right\} $) that has growth
$q=O\left(r^{2}\right)$ as $r\to\infty$. Given a Yang-Mills connection
$A$ with curvature $F$ and structure group $G\subset O\left(N\right)$
($N\geq3$) on the manifold $X$, suppose the self-dual curvature
$F^{+}$ is in the space $L^{2p}$ for some constant $p>1/4$ and
satisfies the pointwise inequality
\begin{equation}
\left|F^{+}\right|\leq\left(1/\gamma\right)\left\{ \left(\left(4p-1\right)/\left(2p^{2}\right)\right)q+\left(1/3\right)S-2\lambda_{\max}\left(W^{+}\right)\right\} \,\,\,\,\,on\,\,\,X\setminus\left\{ x_{0}\right\} .\label{eq:intineq}
\end{equation}
Then either inequality (\ref{eq:intineq}) is an equality (on $X\setminus\left\{ x_{0}\right\} $)
or the self-dual curvature $F^{+}$ is identically zero.
\end{thm}
In the theorem above we also assume that the right hand side of inequality
(\ref{eq:intineq}) is nonnegative and not identically zero. We remark
that the assumption $p>1/4$ is what guarantees that the coefficient
$2-1/\left(2p\right)$ in the proof is positive.
\begin{rem}
\label{rem:singular}When $q$ is smooth on $X$ both theorems above
are valid by changing ``on $X\setminus\left\{ x_{0}\right\} $''
to ``on $X$''.
\end{rem}
\begin{rem}
Both theorems above are valid by changing ``self-dual'', $F^{+}$
and $W^{+}$ to ``anti-self-dual'', $F^{-}$ and $W^{-}$, respectively.
\end{rem}
The two theorems complement each other:
\begin{itemize}
\item In Theorem \ref{thm:vol} we have a volume growth assumption but no
integral assumption.
\item In Theorem \ref{thm:int} we have an integral assumption but no volume
growth assumption.
\end{itemize}
Note that taking $p=1/2$ in Theorem \ref{thm:int} we recover the
pointwise bound of Theorem \ref{thm:vol}. In fact
\[
\frac{4p-1}{2p^{2}}\leq2,
\]
with equality if and only if $p=1/2$. Thus, when $q\geq0$ the pointwise
bound in Theorem \ref{thm:vol} is at least as good as the one in
Theorem \ref{thm:int}. The two theorems may overlap, but neither
set of hypotheses contains the other. For the sphere $S^{4}$ and
the complex projective space $CP^{2}$, both theorems give the same
result. For the Euclidean space $R^{4}$ and the cylinder $S^{3}\times R$,
Theorem \ref{thm:vol} gives a stronger result since Theorem \ref{thm:int}
requires the additional $L^{2p}$ assumption. For the hyperbolic space
$H^{4}$ and the complex hyperbolic space $CH^{2}$, only Theorem
\ref{thm:int} is applicable since these manifolds have exponential
volume growth.

In \cite{BL1981} Bourguignon-Lawson proved an $L^{\infty}$ gap theorem
for $F^{\pm}$ on the sphere $S^{4}$. We improve this result by a
factor $\sqrt{2}$ for the structure group $G=SO\left(3\right)$.
They also proved a $L^{\infty}$ gap theorem for $F^{-}$ on the complex
projective space $CP^{2}$ (Theorem 5.26 in \cite{BL1981}). We improve
this result by a factor $\sqrt{2}$ (resp. $2$) for the structure
group $G\subset SO\left(N\right)$ ($N\geq4$) (resp. $G=SO\left(3\right)$).
In \cite{DM1982} Dodziuk-MinOo proved an $L^{2}$ gap theorem for
$F^{\pm}$ on the Euclidean space $R^{4}$. We prove an $L^{\infty}$
type gap theorem for $F^{\pm}$ on the Euclidean space $R^{4}$ without
any assumption on its integral. In \cite{P1983} Price proved an $L^{2}$
type gap theorem for the full curvature $F$ on the hyperbolic space
$H^{n}$ of dimension $n\geq5$. In \cite{Z2016} Zhou proved an $L^{n/2}$
type gap theorem for the full curvature $F$ on the hyperbolic space
$H^{n}$ of dimension $n\geq16$. We prove an $L^{\infty}$ type gap
theorem for $F^{\pm}$ on the hyperbolic space $H^{4}$. We also prove
an $L^{\infty}$ type gap theorem for $F^{-}$ on the complex hyperbolic
space $CH^{2}$ and an $L^{\infty}$ gap theorem for $F^{\pm}$ on
the cylinder $S^{3}\times R$ without any assumption on its integral.
Applying Theorem \ref{thm:vol} to the manifolds with nonnegative
curvature, applying Theorem \ref{thm:int} to the manifolds with negative
curvature, and using the weighted Poincaré inequalities of Section
2, we get:
\begin{cor}
\label{cor:app}Fix a Lie group $G\subset O\left(N\right)$ ($N\geq3$).
Then:

(i) Given a Yang-Mills connection $A$ with curvature $F$ and structure
group $G$ on the sphere $S^{4}$ (with sectional curvature $1$),
if $F^{\pm}$ satisfies the inequality
\[
\left|F^{\pm}\right|_{L^{\infty}}\leq4/\gamma,
\]
then either $\left|F^{\pm}\right|=4/\gamma$ (everywhere) or $F^{\pm}$
is identically zero.

(ii) Given a Yang-Mills connection $A$ with curvature $F$ and structure
group $G$ on the complex projective space $CP^{2}$ (with holomorphic
sectional curvature $4$), if $F^{-}$ satisfies the inequality
\[
\left|F^{-}\right|_{L^{\infty}}\leq8/\gamma,
\]
then either $\left|F^{-}\right|=8/\gamma$ (everywhere) or $F^{-}$
is identically zero.

(iii) Given a Yang-Mills connection $A$ with curvature $F$ and structure
group $G$ on the cylinder $S^{3}\times R$ (with the product metric),
if $F^{\pm}$ satisfies the inequality
\[
\left|F^{\pm}\right|_{L^{\infty}}\leq2/\gamma,
\]
then either $\left|F^{\pm}\right|=2/\gamma$ (everywhere) or $F^{\pm}$
is identically zero.

(iv) Given a Yang-Mills connection $A$ with curvature $F$ and structure
group $G$ on the Euclidean space $R^{4}$, if $F^{\pm}$ satisfies
the inequality
\[
\left|r^{2}F^{\pm}\right|_{L^{\infty}}\leq2/\gamma,
\]
then $F^{\pm}$ is identically zero.

(v) Given a Yang-Mills connection $A$ with curvature $F$ and structure
group $G$ on the hyperbolic space $H^{4}$ (with sectional curvature
$-1$), if $F^{\pm}$ is in the space $L^{2p}$ for some constant
$3/8\leq p\leq3/4$ and satisfies the pointwise inequality
\[
\left|F^{\pm}\right|\leq\left(1/\gamma\right)\left\{ \left(\left(4p-1\right)/\left(2p^{2}\right)\right)\left(9/4+\left(1/4\right)r^{-2}+\left(3/4\right)\left(\sinh r\right)^{-2}\right)-4\right\} ,
\]
then $F^{\pm}$ is identically zero.

(vi) Given a Yang-Mills connection $A$ with curvature $F$ and structure
group $G$ on the complex hyperbolic space $CH^{2}$ (with holomorphic
sectional curvature $-4$), if $F^{-}$ is in the space $L^{1}$ and
satisfies the pointwise inequality
\[
\left|F^{-}\right|\leq\left(2/\gamma\right)\left(\left(1/4\right)r^{-2}+\left(\sinh r\right)^{-2}-\left(\sinh\left(2r\right)\right)^{-2}\right),
\]
then $F^{-}$ is identically zero.
\end{cor}
To our knowledge, Item (iv) is the first gap theorem on $R^{4}$ that
does not require finite Yang-Mills energy.

Note that Items (i)-(iv) concern nonnegative curvature models and
Items (v)-(vi) concern negatively curved models.

In \cite{BL1981} Bourguignon-Lawson proved that given a non-flat
Yang-Mills connection on the sphere $S^{4}$ with structure group
$G\subset SO\left(4\right)$, if its curvature $F$ satisfies the
inequality $\left|F\right|_{L^{\infty}}\leq4/\gamma$, then the bundle
is one of the four-dimensional bundles of tangent spinors with the
canonical Riemannian connection. Note that $SU\left(2\right)\subset SO\left(4\right)$.
When the structure group is exactly $G=SU\left(2\right)$, we give
an $L^{\infty}$ characterization of the BPST instanton centered at
the origin with unit scale, which appears to be the first characterization
of this instanton in terms of an $L^{\infty}$ extremal property.
\begin{thm}
\label{thm:adhm}Given a non-flat $SU\left(2\right)$ Yang-Mills connection
$A$ on the sphere $S^{4}$ (with sectional curvature $1$), if its
curvature $F$ satisfies the inequality
\[
\left|F\right|_{L^{\infty}\left(S^{4}\right)}\leq4/\gamma,
\]
then, under a stereographic projection, the connection $A$ is gauge
equivalent to the BPST instanton (or the anti-BPST instanton) on the
Euclidean space $R^{4}$ centered at the origin with unit scale.
\end{thm}
In particular, we can characterize the BPST and anti-BPST instantons
centered at the origin with unit scale as the only corresponding $SU\left(2\right)$
instantons on the sphere $S^{4}$ whose curvature $F$ satisfies the
equation $\left|F\right|_{L^{\infty}\left(S^{4}\right)}=4/\gamma$,
up to gauge equivalence. On the other hand, we cannot characterize
these instantons as the only corresponding $SU\left(2\right)$ instantons
on the sphere $S^{4}$ whose curvature $F$ satisfies the equation
$\left|F\right|_{L^{2}\left(S^{4}\right)}=\left(4/\gamma\right)vol\left(S^{4}\right)^{1/2}$.
In fact this equation is satisfied by the curvature of any connection
in the five-parameter family of $SU\left(2\right)$ instantons on
the sphere $S^{4}$ with charge $1$ (see for example the proof of
Theorem \ref{thm:adhm}). In other words, we cannot characterize the
BPST instanton and anti-BPST instanton centered at the origin with
unit scale in terms of the $L^{2}$ norm, but we can characterize
them in terms of the $L^{\infty}$ norm.

The key analytic tool for this work is the combination of a Bochner
inequality (Lemma \ref{lem:bochner}) with the weighted Poincaré inequality.
After multiplication by suitable cutoff functions, we obtain an integral
inequality and conclude that $F^{+}$ is identically zero. The argument
is inspired by vanishing theorems for harmonic forms on complete manifolds
satisfying weighted Poincaré inequalities.

We would like to thank Detang Zhou, Gonçalo Oliveira and Alex Waldron
for their support. We would also like to thank the referee for very
useful suggestions.

\section{Weighted Poincaré inequalities}

\subsection{Model spaces}

First we discuss the weighted Poincaré inequalities that are important
for Corollary \ref{cor:app}.

The weighted Poincaré inequality of the next lemma is known as Hardy's
inequality. Substituting the weight function $q$ of the next lemma
into Theorem \ref{thm:vol}, we get the result for the Euclidean space
$R^{4}$ in Corollary \ref{cor:app}.
\begin{lem}
\label{lem:wpirn}The Euclidean space $R^{n}$ satisfies a weighted
Poincaré inequality with the weight function $q=\left(\left(n-2\right)/2\right)^{2}r^{-2}$,
where $r\left(x\right)=\left|x\right|$.
\end{lem}
The bottom of the spectrum of the Laplacian of the hyperbolic space
$H^{n}$ is $\lambda_{1}\left(\Delta_{H^{n}}\right)=\left(\left(n-1\right)/2\right)^{2}$.
Substituting the weight function $q=\lambda_{1}\left(\Delta_{H^{4}}\right)$
into Theorem \ref{thm:int}, we can get a result for the hyperbolic
space $H^{4}$. Here we use a weight function that improves the bottom
of the spectrum of the Laplacian and leads to a slightly better result
(note that the weight function $q$ of the next lemma satisfies $q>\lambda_{1}\left(\Delta_{H^{n}}\right)$).
In \cite{AK2013} Akutagawa-Kumura proved the weighted Poincaré inequality
of the next lemma, and in \cite{BGG2017} Berchio-Ganguly-Grillo proved
that the inequality is sharp.

Substituting the weight function $q$ of the next lemma into Theorem
\ref{thm:int}, using the fact that the hyperbolic space $H^{4}$
has scalar curvature $S=-12$ and Weyl curvature $W=0$, and using
the fact that for $3/8\leq p\leq3/4$ we have $\left(\left(4p-1\right)/\left(2p^{2}\right)\right)q-4>0$,
we get the result for the hyperbolic space $H^{4}$ in Corollary \ref{cor:app}.
\begin{lem}
\label{lem:wpihn}The hyperbolic space $H^{n}$ (with sectional curvature
$-1$) satisfies a weighted Poincaré inequality with the weight function
\[
q=\left(\left(n-1\right)/2\right)^{2}+\left(1/4\right)r^{-2}+\left\{ \left(n-1\right)\left(n-3\right)/4\right\} \left(\sinh r\right)^{-2},
\]
where $r\left(\cdot\right)=dist_{H^{n}}\left(\cdot,x_{0}\right)$.
\end{lem}
The bottom of the spectrum of the Laplacian of the complex hyperbolic
space $CH^{n}$ is $\lambda_{1}\left(\Delta_{CH^{n}}\right)=n^{2}$
(see for example Li-Wang \cite{LW2005}). But substituting the weight
function $q=\lambda_{1}\left(\Delta_{CH^{2}}\right)$ into Theorem
\ref{thm:int}, we do not obtain a nontrivial result for the complex
hyperbolic space $CH^{2}$. Indeed taking the optimal choice $p=1/2$
in Theorem \ref{thm:int} and using the fact that $S=-24$ and $W^{-}=0$
the pointwise condition becomes
\[
\left|F^{-}\right|\leq\frac{1}{\gamma}\left(2q-8\right),
\]
so the theorem gives no nontrivial gap result for $q=\lambda_{1}\left(\Delta_{CH^{2}}\right)=4$.
To overcome this difficulty we find a weighted Poincaré inequality
that improves the bottom of the spectrum of the Laplacian, following
closely the ideas in \cite{AK2013}. Note that the weight function
$q$ of the next lemma satisfies $q>\lambda_{1}\left(\Delta_{CH^{n}}\right)$,
and in particular $2q-8>0$ for $CH^{2}$. The weighted Poincaré inequality
in the next lemma is the complex hyperbolic analogue of the weighted
Poincaré inequality in Lemma \ref{lem:wpihn}, since both are obtained
from the same general construction of Akutagawa-Kumura \cite{AK2013}.
The inequality below appeared in an earlier version of the present
work. Subsequently, the work Fischer-Peyerimhoff \cite{FP2025} implies
the sharpness of the corresponding inequality as a consequence of
their optimal Hardy inequality for Damek-Ricci spaces. Indeed, up
to rescaling the metric, $CH^{n}$ is the Damek-Ricci space $X^{2\left(n-1\right),1}$,
and their main theorem gives precisely the weight in the next lemma.

Substituting the weight function $q$ of the next lemma and $p=1/2$
into Theorem \ref{thm:int} and using the fact that the complex hyperbolic
space $CH^{2}$ has scalar curvature $S=-24$ and Weyl curvature $W^{-}=0$,
we get the result for the complex hyperbolic space $CH^{2}$ in Corollary
\ref{cor:app}.
\begin{lem}
\label{lem:wpichn}The complex hyperbolic space $CH^{n}$ (with holomorphic
sectional curvature $-4$) satisfies a weighted Poincaré inequality
with the weight function
\[
q=n^{2}+\left(1/4\right)r^{-2}+\left(n-1\right)^{2}\left(\sinh r\right)^{-2}-\left(\sinh\left(2r\right)\right)^{-2},
\]
where $r\left(\cdot\right)=dist_{CH^{n}}\left(\cdot,x_{0}\right)$.
\end{lem}
\begin{proof}
In \cite{AK2013} Akutagawa-Kumura proved that any complete Riemannian
manifold with a pole $x_{0}$ satisfies a weighted Poincaré inequality
with the weight function
\[
q=\left(1/4\right)r^{-2}+\left(1/4\right)\left(\Delta r\right)^{2}-\left(1/2\right)\left|\nabla\nabla r\right|^{2}-\left(1/2\right)Ric\left(\nabla r,\nabla r\right),
\]
where $r\left(\cdot\right)=dist\left(\cdot,x_{0}\right)$. Substituting
the distance function $r$ into Bochner's formula and using the fact
that $\left|\nabla r\right|^{2}=1$, we see that the sum of the last
two terms on the right hand side is equal to $\left(1/2\right)\left\langle \nabla\Delta r,\nabla r\right\rangle $,
so
\[
q=\left(1/4\right)r^{-2}+\left(1/4\right)\left(\Delta r\right)^{2}+\left(1/2\right)\left\langle \nabla\Delta r,\nabla r\right\rangle .
\]
In the special case of the complex hyperbolic space $CH^{n}$ the
Laplacian of the distance function $r$ is (see for example Theorem
1.6 in Li-Wang \cite{LW2005})
\[
\Delta r=2\left(n-1\right)\coth r+2\coth\left(2r\right),
\]
so
\[
\left(1/4\right)\left(\Delta r\right)^{2}=n^{2}+n\left(n-1\right)\left(\sinh r\right)^{-2}+\left(\sinh\left(2r\right)\right)^{-2},
\]
\[
\left(1/2\right)\left\langle \nabla\Delta r,\nabla r\right\rangle =\left(1/2\right)\frac{d}{dr}\left(\Delta r\right)=-\left(n-1\right)\left(\sinh r\right)^{-2}-2\left(\sinh\left(2r\right)\right)^{-2}.
\]
Substituting these two terms into the equation of the weight function
$q$, we get the result.
\end{proof}

\subsection{Other geometric settings}

Next we give more examples of weighted Poincaré inequalities. Using
these inequalities, we could get more applications of Theorem \ref{thm:vol}
and Theorem \ref{thm:int}.

Manifolds with positive spectrum satisfy a weighted Poincaré inequality.
\begin{example}
The bottom of the spectrum of the Laplacian of a Riemannian manifold
$X$ is defined by $\lambda_{1}\left(\Delta\right)=\inf_{\phi}\int_{X}\left|\nabla\phi\right|^{2}$,
where the infimum is taken over all smooth functions $\phi$ on the
manifold $X$ with compact support and $\int_{X}\phi^{2}=1$. In particular
the manifold $X$ satisfies a weighted Poincaré inequality with the
weight function $q=\lambda_{1}\left(\Delta\right)$.
\end{example}
Stable minimal hypersurfaces in a general Riemannian manifold satisfy
a weighted Poincaré inequality.
\begin{example}
A minimal hypersurface $X^{n}$ in a Riemannian manifold $Y^{n+1}$
is called stable if the second variation of the area functional is
nonnegative for any normal variation with compact support. It is well
known that the hypersurface $X^{n}$ is stable if and only if it satisfies
a weighted Poincaré inequality with the weight function $q=\left|A\right|^{2}+Ric_{Y}\left(N,N\right)$.
Here we have the second fundamental form $A$ of the hypersurface
and the Ricci curvature $Ric_{Y}$ of the ambient space applied to
the unit normal vector $N$.
\end{example}
Minimal hypersurfaces in the Euclidean space satisfy a weighted Poincaré
inequality.
\begin{example}
In \cite{C1997} Carron proved that a minimal hypersurface $X^{n}$
in the Euclidean space $R^{n+1}$ satisfies a weighted Poincaré inequality
with the weight function $q=\left(\left(n-2\right)/2\right)^{2}r^{-2}$,
where $r\left(\cdot\right)=dist_{R^{n+1}}\left(\cdot,x_{0}\right)|X$.
\end{example}
Finally we give more references on weighted Poincaré inequalities.
For applications of weighted Poincaré inequalities to the geometry
and topology of Riemannian manifolds and submanifolds see for example
Carron \cite{C1997}, Li-Wang \cite{LW2006}, Cheng-Zhou \cite{CZ2009},
Lam \cite{L2010} and Vieira \cite{V2016}. In \cite{H2011}, Corollary
1.3 with $\alpha=1$, Hein proved a weighted Poincaré inequality in
certain asymptotically conical manifolds that may lead to new applications
of Theorem \ref{thm:vol}. We would like to find weighted Poincaré
inequalities in interesting ALE manifolds (for example the Eguchi-Hanson
manifold) and ALF manifolds (for example the Euclidean Schwarzschild
manifold) leading to new applications of Theorem \ref{thm:vol}.

\section{Proofs}

\subsection{Bochner formula}

In the proof of Theorem \ref{thm:vol} and Theorem \ref{thm:int}
we use a Bochner type formula based on Bourguignon-Lawson \cite{BL1981}
and Gursky-Kelleher-Streets \cite{GKS2018} (for details see \cite{V2024}).
Although the constant in the term of the Weyl curvature is different,
the result follows easily from the proof in \cite{V2024}.
\begin{lem}
\label{lem:bochner}Suppose we have a four-dimensional Riemannian
manifold $X$ with scalar curvature $S$ and Weyl curvature $W$,
and a Yang-Mills connection $A$ with curvature $F$ and structure
group $G\subset O\left(N\right)$ ($N\geq3$) on the manifold $X$.
Fix a constant $p>0$. Then the self-dual curvature $F^{+}$ satisfies
the pointwise inequality
\begin{align*}
\left|F^{+}\right|^{p}\Delta\left|F^{+}\right|^{p}\geq\left(1-1/\left(2p\right)\right)\left|\nabla\left|F^{+}\right|^{p}\right|^{2}+\left(p/3\right)S\left|F^{+}\right|^{2p}\\
-2p\lambda_{\max}\left(W^{+}\right)\left|F^{+}\right|^{2p}-p\gamma\left|F^{+}\right|^{2p+1}.
\end{align*}
Here we use the inner products
\[
\left|F^{+}\right|^{2}=\sum_{1\leq i<j\leq4}\left|F_{ij}^{+}\right|^{2},\,\,\,\,\,\left\langle A,B\right\rangle _{so\left(N\right)}=-c\cdot tr\left(AB\right),\,\,\,\,\,c>0,
\]
and we denote
\[
\gamma=\begin{cases}
4/\sqrt{12c}, & N=3,\\
4/\sqrt{6c}, & N\geq4.
\end{cases}
\]
Recall that $F^{+}=1/2\left(F+*F\right)$ and $\lambda_{\max}\left(W^{+}\right)$
is the largest eigenvalue of $W^{+}=1/2\left(W+*W\right)$.
\end{lem}

\subsection{Proof of Theorem \ref{thm:vol}}

First we combine an estimate obtained from Lemma \ref{lem:bochner}
with an estimate obtained from the weighted Poincaré inequality. By
Lemma \ref{lem:bochner}, we have
\[
-\left|F^{+}\right|^{1/2}\Delta\left|F^{+}\right|^{1/2}+\left(1/6\right)S\left|F^{+}\right|-\lambda_{\max}\left(W^{+}\right)\left|F^{+}\right|-\left(\gamma/2\right)\left|F^{+}\right|^{2}\leq0.
\]
Recall that $r\left(\cdot\right)=dist\left(\cdot,x_{0}\right)$ and
$B_{R}=\left\{ r<R\right\} $. Take the logarithmic cutoff function
$\phi$ on the manifold given by
\[
\phi=\begin{cases}
1 & B_{R},\\
2-\log r/\log R & B_{R^{2}}\setminus B_{R},\\
0 & X\setminus B_{R^{2}}.
\end{cases}
\]
This cutoff function is often used in parabolicity arguments and minimal
surface theory and it is important in the next part of the proof.
Multiplying this inequality by the function $\phi^{2}$ and integrating
by parts, we get
\begin{align*}
 & \int\left|\nabla\left|F^{+}\right|^{1/2}\right|^{2}\phi^{2}+2\int\left|F^{+}\right|^{1/2}\phi\left\langle \nabla\left|F^{+}\right|^{1/2},\nabla\phi\right\rangle +\left(1/6\right)\int S\left|F^{+}\right|\phi^{2}\\
 & -\int\lambda_{\max}\left(W^{+}\right)\left|F^{+}\right|\phi^{2}-\left(\gamma/2\right)\int\left|F^{+}\right|^{2}\phi^{2}\leq0.
\end{align*}
On the other hand, substituting the compactly supported function $\left|F^{+}\right|^{1/2}\phi$
into the weighted Poincaré inequality and rearranging the terms, we
have
\begin{align*}
 & -\int\left|\nabla\left|F^{+}\right|^{1/2}\right|^{2}\phi^{2}-2\int\left|F^{+}\right|^{1/2}\phi\left\langle \nabla\left|F^{+}\right|^{1/2},\nabla\phi\right\rangle \\
 & +\int q\left|F^{+}\right|\phi^{2}\leq\int\left|F^{+}\right|\left|\nabla\phi\right|^{2}.
\end{align*}
Summing the above estimates, we get
\[
\int\left\{ q+\left(1/6\right)S-\lambda_{\max}\left(W^{+}\right)-\left(\gamma/2\right)\left|F^{+}\right|\right\} \left|F^{+}\right|\phi^{2}\leq\int\left|F^{+}\right|\left|\nabla\phi\right|^{2}.
\]

Next we estimate the right hand side of this inequality. We denote
by $C$ a positive constant independent of $R$, which may change
from line to line. Substituting the gradient of the function $\phi$
(which vanishes outside the ring $B_{R^{2}}\setminus B_{R}$), writing
the ring $B_{R^{2}}\setminus B_{R}$ as the disjoint union $B_{R^{2}}\setminus B_{R}=\bigcup_{i=\log R}^{2\log R-1}B_{e^{i+1}}\setminus B_{e^{i}}$,
and using the fact that the self-dual curvature $F^{+}$ has growth
$\left|F^{+}\right|=O\left(r^{2-\alpha}\right)$, we see that for
$R$ sufficiently large
\begin{align*}
 & \int\left\{ q+\left(1/6\right)S-\lambda_{\max}\left(W^{+}\right)-\left(\gamma/2\right)\left|F^{+}\right|\right\} \left|F^{+}\right|\phi^{2}\\
 & \leq C\left(\log R\right)^{-2}\sum_{i=\log R}^{2\log R-1}\int_{B_{e^{i+1}}\setminus B_{e^{i}}}r^{-\alpha}.
\end{align*}
Since the manifold has volume growth $vol\left(B_{R}\right)=O\left(R^{\alpha}\right)$,
we have $\int_{B_{e^{i+1}}\setminus B_{e^{i}}}r^{-\alpha}=O\left(1\right)$,
so for $R$ sufficiently large
\[
\int\left\{ q+\left(1/6\right)S-\lambda_{\max}\left(W^{+}\right)-\left(\gamma/2\right)\left|F^{+}\right|\right\} \left|F^{+}\right|\phi^{2}\leq C\left(\log R\right)^{-1}.
\]

Finally we show that either inequality (\ref{eq:volineq}) is an equality
(on $X\setminus\left\{ x_{0}\right\} $) or the self-dual curvature
$F^{+}$ is identically zero. By assumption, the term inside the braces
is nonnegative on $X\setminus\left\{ x_{0}\right\} $. Taking the
limit as $R\to\infty$, we get
\[
\left\{ q+\left(1/6\right)S-\lambda_{\max}\left(W^{+}\right)-\left(\gamma/2\right)\left|F^{+}\right|\right\} \left|F^{+}\right|=0\,\,\,\,\,on\,\,\,X\setminus\left\{ x_{0}\right\} .
\]
Suppose that inequality (\ref{eq:volineq}) is not an equality (on
$X\setminus\left\{ x_{0}\right\} $). In this case the term inside
the braces is positive at some point of $X\setminus\left\{ x_{0}\right\} $,
so by continuity the self-dual curvature $F^{+}$ vanishes in an open
subset of the manifold. By the unique continuation principle, we conclude
that the self-dual curvature $F^{+}$ is identically zero.

\subsection{Proof of Theorem \ref{thm:int}}

First, as in the proof of Theorem \ref{thm:vol}, we combine an estimate
obtained from Lemma \ref{lem:bochner} with an estimate obtained from
the weighted Poincaré inequality, but the calculation here is more
involved. By Lemma \ref{lem:bochner}, we have
\begin{align*}
 & -\left|F^{+}\right|^{p}\Delta\left|F^{+}\right|^{p}+\left(1-1/\left(2p\right)\right)\left|\nabla\left|F^{+}\right|^{p}\right|^{2}+\left(p/3\right)S\left|F^{+}\right|^{2p}\\
 & -2p\lambda_{\max}\left(W^{+}\right)\left|F^{+}\right|^{2p}-p\gamma\left|F^{+}\right|^{2p+1}\leq0.
\end{align*}
Recall that $r\left(\cdot\right)=dist\left(\cdot,x_{0}\right)$ and
$B_{R}=\left\{ r<R\right\} $. Take the cutoff function $\phi$ on
the manifold given by
\[
\phi=\begin{cases}
1 & B_{R},\\
2-r/R & B_{2R}\setminus B_{R},\\
0 & X\setminus B_{2R}.
\end{cases}
\]
Multiplying this inequality by the function $\phi^{2}$, integrating
by parts and writing $\tilde{p}=2-1/\left(2p\right)$ (which is positive
since $p>1/4$), we get
\begin{align*}
 & \tilde{p}\int\left|\nabla\left|F^{+}\right|^{p}\right|^{2}\phi^{2}+2\int\left|F^{+}\right|^{p}\phi\left\langle \nabla\left|F^{+}\right|^{p},\nabla\phi\right\rangle +\left(p/3\right)\int S\left|F^{+}\right|^{2p}\phi^{2}\\
 & -2p\int\lambda_{\max}\left(W^{+}\right)\left|F^{+}\right|^{2p}\phi^{2}-p\gamma\int\left|F^{+}\right|^{2p+1}\phi^{2}\leq0.
\end{align*}
Fix a small $\delta>0$. Applying the inequality $2ab\leq\delta a^{2}+\delta^{-1}b^{2}$
to the second term on the left hand side and rearranging terms, we
get (i):
\begin{align*}
 & \left(\tilde{p}-\delta\right)\int\left|\nabla\left|F^{+}\right|^{p}\right|^{2}\phi^{2}+\left(p/3\right)\int S\left|F^{+}\right|^{2p}\phi^{2}-2p\int\lambda_{\max}\left(W^{+}\right)\left|F^{+}\right|^{2p}\phi^{2}\\
 & -p\gamma\int\left|F^{+}\right|^{2p+1}\phi^{2}\leq\delta^{-1}\int\left|F^{+}\right|^{2p}\left|\nabla\phi\right|^{2}.
\end{align*}
On the other hand, substituting the compactly supported function $\left|F^{+}\right|^{p}\phi$
into the weighted Poincaré inequality and using the inequality $2ab\leq\delta a^{2}+\delta^{-1}b^{2}$
as before, we have
\[
\int q\left|F^{+}\right|^{2p}\phi^{2}\leq\left(1+\delta\right)\int\left|\nabla\left|F^{+}\right|^{p}\right|^{2}\phi^{2}+\left(1+\delta^{-1}\right)\int\left|F^{+}\right|^{2p}\left|\nabla\phi\right|^{2}.
\]
We denote by $C$ a positive constant independent of $R$ and $\delta$,
which may change from line to line. Multiplying this inequality by
$\left(\tilde{p}-\delta\right)\left(1+\delta\right)^{-1}$, using
the inequalities $\left(\tilde{p}-\delta\right)\left(1+\delta\right)^{-1}\geq\tilde{p}-C\delta$
and $\left(1+\delta^{-1}\right)\left(\tilde{p}-\delta\right)\left(1+\delta\right)^{-1}\leq C\delta^{-1}$,
and rearranging terms, we get (ii):
\begin{align*}
 & -\left(\tilde{p}-\delta\right)\int\left|\nabla\left|F^{+}\right|^{p}\right|^{2}\phi^{2}+\tilde{p}\int q\left|F^{+}\right|^{2p}\phi^{2}\\
 & \leq C\left(\delta\int q\left|F^{+}\right|^{2p}\phi^{2}+\delta^{-1}\int\left|F^{+}\right|^{2p}\left|\nabla\phi\right|^{2}\right).
\end{align*}
Summing the above estimates ((i) and (ii)), we get
\begin{align*}
 & \int\left\{ \left(2-1/\left(2p\right)\right)q+\left(p/3\right)S-2p\lambda_{\max}\left(W^{+}\right)-p\gamma\left|F^{+}\right|\right\} \left|F^{+}\right|^{2p}\phi^{2}\\
 & \leq C\left(\delta\int q\left|F^{+}\right|^{2p}\phi^{2}+\delta^{-1}\int\left|F^{+}\right|^{2p}\left|\nabla\phi\right|^{2}\right).
\end{align*}

Next, as in the proof of Theorem \ref{thm:vol}, we estimate the right
hand side of this inequality. The choice of $\delta$ is tricky. Fix
a small $\epsilon>0$ and take $\delta=\epsilon R^{-2}$. Substituting
the function $\phi$ (which vanishes outside the ball $B_{2R}$) and
its gradient (which vanishes outside the ring $B_{2R}\setminus B_{R}$),
and using the fact that the weight function $q$ has growth $q=O\left(r^{2}\right)$
(say $q\leq Cr^{2}$ for $r\geq R_{0}$), we get
\begin{align*}
 & \int\left\{ \left(2-1/\left(2p\right)\right)q+\left(p/3\right)S-2p\lambda_{\max}\left(W^{+}\right)-p\gamma\left|F^{+}\right|\right\} \left|F^{+}\right|^{2p}\phi^{2}\\
 & \leq C\left(\epsilon R^{-2}\int_{B_{R_{0}}}q\left|F^{+}\right|^{2p}+\epsilon\int_{B_{2R}\setminus B_{R_{0}}}\left|F^{+}\right|^{2p}+\epsilon^{-1}\int_{B_{2R}\setminus B_{R}}\left|F^{+}\right|^{2p}\right).
\end{align*}

Finally, as in the proof of Theorem \ref{thm:vol}, we show that either
inequality (\ref{eq:intineq}) is an equality (on $X\setminus\left\{ x_{0}\right\} $)
or the self-dual curvature $F^{+}$ is identically zero. By assumption,
the term inside the braces on the left hand side is nonnegative on
$X\setminus\left\{ x_{0}\right\} $ and the self-dual curvature $F^{+}$
is in the space $L^{2p}$. Taking the limit as $R\to\infty$ and then
taking the limit as $\epsilon\to0$, we see that the right hand side
becomes zero, so
\[
\left\{ \left(2-1/\left(2p\right)\right)q+\left(p/3\right)S-2p\lambda_{\max}\left(W^{+}\right)-p\gamma\left|F^{+}\right|\right\} \left|F^{+}\right|^{2p}=0\,\,\,\,\,on\,\,\,X\setminus\left\{ x_{0}\right\} .
\]
As in the proof of Theorem \ref{thm:vol}, if inequality (\ref{eq:intineq})
is not an equality (on $X\setminus\left\{ x_{0}\right\} $), then
the self-dual curvature $F^{+}$ is identically zero.

\subsection{Proof of Corollary \ref{cor:app}}

We proved the results for Euclidean space $R^{4}$, the hyperbolic
space $H^{4}$ and the complex hyperbolic space $CH^{2}$ above Lemma
\ref{lem:wpirn}, Lemma \ref{lem:wpihn} and Lemma \ref{lem:wpichn},
respectively. Substituting the weight function $q=0$ into Theorem
\ref{thm:vol} (recalling Remark \ref{rem:singular}) and using the
fact that the sphere $S^{4}$ has scalar curvature $S=12$ and Weyl
curvature $W=0$, the complex projective space $CP^{2}$ has scalar
curvature $S=24$ and Weyl curvature $W^{-}=0$, and the cylinder
$S^{3}\times R$ has scalar curvature $S=6$ and Weyl curvature $W=0$,
we get the results for the sphere $S^{4}$, the complex projective
space $CP^{2}$ and the cylinder $S^{3}\times R$.

\subsection{Proof of Theorem \ref{thm:adhm}}

Note that $SU\left(2\right)\subset SO\left(4\right)$. Taking $c=1/2$
(without loss of generality) and $N=4$ in the definition of $\gamma$
(see Lemma \ref{lem:bochner}), we see that the inequality in the
assumption of the theorem becomes the inequality $\left|F\right|_{L^{\infty}\left(S^{4}\right)}\leq\sqrt{3}$.

First we show that the connection $A$ is an instanton, the curvature
$F$ satisfies the equation $\left|F\right|_{g_{S^{4}}}=\sqrt{3}$
(everywhere) and the absolute value of the charge is $1$. The proof
is as follows. Since $\left|F\right|^{2}=\left|F^{+}\right|^{2}+\left|F^{-}\right|^{2}$,
we get $\left|F^{+}\right|_{L^{\infty}\left(S^{4}\right)}\leq\sqrt{3}$
and $\left|F^{-}\right|_{L^{\infty}\left(S^{4}\right)}\leq\sqrt{3}$.
Applying Corollary \ref{cor:app} to the self-dual curvature $F^{+}$
(resp. anti-self-dual curvature $F^{-}$), we find that $\left|F^{+}\right|_{g_{S^{4}}}=\sqrt{3}$
everywhere or $F^{+}$ is identically zero (resp. $\left|F^{-}\right|_{g_{S^{4}}}=\sqrt{3}$
everywhere or $F^{-}$ is identically zero). Since $\left|F\right|^{2}=\left|F^{+}\right|^{2}+\left|F^{-}\right|^{2}$
and $\left|F\right|_{g_{S^{4}}}\leq\sqrt{3}$, the equalities $\left|F^{+}\right|_{g_{S^{4}}}=\sqrt{3}$
and $\left|F^{-}\right|_{g_{S^{4}}}=\sqrt{3}$ cannot hold simultaneously.
Since $A$ is non-flat, we find that either (i) $\left|F^{+}\right|_{g_{S^{4}}}=\sqrt{3}$
everywhere and $F^{-}$ is identically zero or (ii) $F^{+}$ is identically
zero and $\left|F^{-}\right|_{g_{S^{4}}}=\sqrt{3}$ everywhere. In
both cases we conclude that the connection $A$ is an instanton and
the curvature $F$ satisfies $\left|F\right|_{g_{S^{4}}}=\sqrt{3}$
everywhere. Since $vol\left(S^{4}\right)=8\pi^{2}/3$, we get that
the absolute value of the charge is $\left(8\pi^{2}\right)^{-1}\left|F\right|_{L^{2}\left(S^{4}\right)}^{2}=1$.

Next we discuss the corresponding instanton on the Euclidean space
$R^{4}$. We consider the case $F^{-}=0$; the case $F^{+}=0$ gives
the anti-BPST instanton centered at the origin with unit scale by
the same argument. Using the fact that the connection $A$ is an $SU\left(2\right)$
instanton on the sphere $S^{4}$ with charge $1$, we see that under
the stereographic projection of $S^{4}\setminus\left\{ N\right\} $
the connection $A$ corresponds to an $SU\left(2\right)$ instanton
on the Euclidean space $R^{4}$ with charge $1$. Using the fact that
the curvature $F$ satisfies the equation $\left|F\right|_{g_{S^{4}}}=\sqrt{3}$
everywhere and writing the metric of the sphere $S^{4}$ in terms
of the stereographic projection as $g_{S^{4}}=4\left(1+\left|x\right|^{2}\right)^{-2}g_{R^{4}}$,
we see that the curvature of the corresponding instanton (also denoted
by $F$) satisfies the equation
\[
\left|F\right|_{g_{R^{4}}}^{2}\left(x\right)=48\left(1+\left|x\right|^{2}\right)^{-4}.
\]

Finally we show that the corresponding instanton is the BPST instanton
centered at the origin with unit scale. It is well known that, up
to gauge equivalence, $SU\left(2\right)$ instantons on the Euclidean
space $R^{4}$ with charge $1$ are an explicit five-parameter family
of connections $\left\{ A_{x_{0},a}\right\} $ determined by a center
$x_{0}$ in $R^{4}$ and a scale $a>0$, and the curvature $F\left(A_{x_{0},a}\right)$
satisfies the equation
\[
\left|F\left(A_{x_{0},a}\right)\right|_{g_{R^{4}}}^{2}\left(x\right)=48a^{4}\left(a^{2}+\left|x-x_{0}\right|^{2}\right)^{-4}.
\]
See for example Chapter 3 in \cite{DK1990} and Section 2 in \cite{VN2008}.
Comparing the above equations, we conclude that the connection $A$
is gauge equivalent to the instanton $A_{x_{0},a}$ with center $x_{0}=0$
and scale $a=1$.
\begin{rem}
The first part of the proof of Theorem \ref{thm:adhm} can also be
obtained from the $L^{2}$ gap theorem of Gursky-Kelleher-Streets
\cite{GKS2018}. The $L^{\infty}$ bound implies $\left|F\right|_{L^{2}\left(S^{4}\right)}\leq8\pi^{2}$
(which becomes $16\pi^{2}$ in the normalization of \cite{GKS2018}).
Then the $L^{2}$ gap theorem in \cite{GKS2018} implies that $A$
is an instanton. We have chosen to give the proof above since it follows
directly from the $L^{\infty}$ gap theorem developed in the present
paper.
\end{rem}

\lyxaddress{Departamento de Matemática, Universidade Federal do Espírito Santo,
Vitória, ES, Brazil. Email: matheus.vieira@ufes.br}
\end{document}